\documentclass[12pt]{article}
\newcommand{\be}{\begin{equation}}
\newcommand{\ee}{\end{equation}}
\newcommand{\bea}{\begin{eqnarray}}
\newcommand{\eea}{\end{eqnarray}}
\newcommand{\ba}{\begin{array}}
\newcommand{\ea}{\end{array}}

\newcommand{\bc}{\begin{center}}
\newcommand{\ec}{\end{center}}
\newcommand{\ben}{\begin{enumerate}}
\newcommand{\een}{\end{enumerate}}
\newcommand{\bfi}{\begin{figure}}
\newcommand{\efi}{\end{figure}}

\newcommand{\bq}{\begin{quote}}
\newcommand{\eq}{\end{quote}}
\newcommand{\bqu}{\begin{quotation}}
\newcommand{\equ}{\end{quotation}}
\newenvironment{emphit}{\begin{itemize}}{\end{itemize}}
\newcommand{\bemp}{\begin{emphit}}
\newcommand{\eemp}{\end{emphit}}

\newcommand{\bt}{\begin{tabular}}
\newcommand{\et}{\end{tabular}}

\newtheorem{myth}{Theorem}[section]
\newtheorem{mylem}{Lemma}[section]
\newtheorem{mycor}{Corollary}[section]

\newtheorem{mydef}{Definition}[section]

\begin{document}
\date{}
\title{ON THE UNIVERSALITY OF SOME SMARANDACHE LOOPS OF BOL-MOUFANG TYPE
}
\author{T\`em\'it\'op\'e Gb\'ol\'ah\`an Ja\'iy\'e\d ol\'a\thanks{On Doctorate Programme at
the University of Agriculture Abeokuta, Nigeria.}
\thanks{All correspondence to be addressed to this author}\\
Department of Mathematics,\\
Obafemi Awolowo University, Ile Ife, Nigeria.\\
jaiyeolatemitope@yahoo.com, tjayeola@oauife.edu.ng} \maketitle

\begin{abstract}
A Smarandache quasigroup(loop) is shown to be universal if all its
$f,g$-principal isotopes are Smarandache $f,g$-principal isotopes.
Also, weak Smarandache loops of Bol-Moufang type such as
Smarandache: left(right) Bol, Moufang and extra loops are shown to
be universal if all their $f,g$-principal isotopes are Smarandache
$f,g$-principal isotopes. Conversely, it is shown that if these weak
Smarandache loops of Bol-Moufang type are universal, then some
autotopisms are true in the weak Smarandache sub-loops of the weak
Smarandache loops of Bol-Moufang type relative to some Smarandache
elements. Futhermore, a Smarandache left(right) inverse property
loop in which all its $f,g$-principal isotopes are Smarandache
$f,g$-principal isotopes is shown to be universal if and only if it
is a Smarandache left(right) Bol loop in which all its
$f,g$-principal isotopes are Smarandache $f,g$-principal isotopes.
Also, it is established that a Smarandache inverse property loop in
which all its $f,g$-principal isotopes are Smarandache
$f,g$-principal isotopes is universal if and only if it is a
Smarandache Moufang loop in which all its $f,g$-principal isotopes
are Smarandache $f,g$-principal isotopes. Hence, some of the
autotopisms earlier mentioned are found to be true in the
Smarandache sub-loops of universal Smarandache: left(right) inverse
property loops and inverse property loops.
\end{abstract}

\section{Introduction}
W. B. Vasantha Kandasamy initiated the study of Smarandache loops
(S-loop) in 2002. In her book \cite{phd75}, she defined a
Smarandache loop (S-loop) as a loop with at least a subloop which
forms a subgroup under the binary operation of the loop called a
Smarandache subloop (S-subloop). In \cite{sma2}, the present author
defined a Smarandache quasigroup (S-quasigroup) to be a quasigroup
with at least a non-trivial associative subquasigroup called a
Smarandache subquasigroup (S-subquasigroup). Examples of Smarandache
quasigroups are given in Muktibodh \cite{muk}. For more on
quasigroups, loops and their properties, readers should check
\cite{phd3}, \cite{phd41},\cite{phd39}, \cite{phd49}, \cite{phd42}
and \cite{phd75}. In her (W.B. Vasantha Kandasamy) first paper
\cite{phd83}, she introduced Smarandache : left(right) alternative
loops, Bol loops, Moufang loops, and Bruck loops. But in
\cite{sma1}, the present author introduced Smarandache : inverse
property loops (IPL), weak inverse property loops (WIPL), G-loops,
conjugacy closed loops (CC-loop), central loops, extra loops,
A-loops, K-loops, Bruck loops, Kikkawa loops, Burn loops and
homogeneous loops. The isotopic invariance of types and varieties of
quasigroups and loops described by one or more equivalent
identities, especially those that fall in the class of Bol-Moufang
type loops as first named by Fenyves \cite{phd56} and \cite{phd50}
in the 1960s  and later on in this $21^{st}$ century by Phillips and
Vojt\v echovsk\'y \cite{phd9}, \cite{phd61} and \cite{phd124} have
been of interest to researchers in loop theory in the recent past.
For example, loops such as Bol loops, Moufang loops, central loops
and extra loops are the most popular loops of Bol-Moufang type whose
isotopic invariance have been considered. Their identities relative
to quasigroups and loops have also been investigated by Kunen
\cite{ken1} and \cite{ken2}. A loop is said to be universal relative
to a property ${\cal P}$ if it is isotopic invariant relative to
${\cal P}$, hence such a loop is called a universal ${\cal P}$ loop.
This language is well used in \cite{phd88}. The universality of most
loops of Bol-Moufang types have been studied as summarised in
\cite{phd3}. Left(Right) Bol loops, Moufang loops, and extra loops
have all been found to be isotopic invariant. But some types of
central loops were shown to be universal in Ja\'iy\'e\d ol\'a
\cite{tope} and \cite{phdtope} under some conditions. Some other
types of loops such as A-loops, weak inverse property loops and
cross inverse property loops (CIPL) have been found be universal
under some neccessary and sufficient conditions in \cite{phd40},
\cite{phd43} and \cite{phd30} respectively. Recently, Michael Kinyon
et. al. \cite{phd95}, \cite{phd118}, \cite{phd119} solved the
Belousov problem concerning the universality of F-quasigroups which
has been open since 1967 by showing that all the isotopes of
F-quasigroups are Moufang loops.

In this work, the universality of the Smarandache concept in loops
is investigated. That is, will all isotopes of an S-loop be an
S-loop? The answer to this could be 'yes' since every isotope of a
group is a group (groups are G-loops). Also, the universality of
weak Smarandache loops, such as Smarandache Bol loops (SBL),
Smarandache Moufang loops (SML) and Smarandache extra loops (SEL)
will also be investigated despite the fact that it could be expected
to be true since Bol loops, Moufang loops and extra loops are
universal. The universality of a Smarandache inverse property loop
(SIPL) will also be considered.

\section{Preliminaries}
\begin{mydef}
A loop is called a Smarandache left inverse property loop (SLIPL) if
it has at least a non-trivial subloop with the LIP.

A loop is called a Smarandache right inverse property loop (SRIPL)
if it has at least a non-trivial subloop with the RIP.

A loop is called a Smarandache inverse property loop (SIPL) if it
has at least a non-trivial subloop with the IP.

A loop is called a Smarandache right Bol-loop (SRBL) if it has at
least a non-trivial subloop that is a right Bol(RB)-loop.

A loop is called a Smarandache left Bol-loop (SLBL) if it has at
least a non-trivial subloop that is a left Bol(LB)-loop.

A loop is called a Smarandache central-loop (SCL) if it has at least
a non-trivial subloop that is a central-loop.

A loop is called a Smarandache extra-loop (SEL) if it has at least a
non-trivial subloop that is a extra-loop.

A loop is called a Smarandache A-loop (SAL) if it has at least a
non-trivial subloop that is a A-loop.

A loop is called a Smarandache Moufang-loop (SML) if it has at least
a non-trivial subloop that is a Moufang-loop.
\end{mydef}

\begin{mydef}
Let $(G,\oplus)$ and $(H,\otimes)$ be two distinct quasigroups. The
triple $(A,B,C)$ such that $A,B,C~:~(G,\oplus)\rightarrow
(H,\otimes)$ are bijections is said to be an isotopism if and only
if
\begin{displaymath}
xA\otimes yB=(x\oplus y)C~\forall~x,y\in G.
\end{displaymath}
Thus, $H$ is called an isotope of $G$ and they are said to be
isotopic. If $C=I$, then the triple is called a principal isotopism
and $(H,\otimes)=(G,\otimes )$ is called a principal isotope of
$(G,\oplus )$. If in addition, $A=R_g$, $B=L_f$, then the triple is
called an $f,g$-principal isotopism, thus $(G,\otimes )$ is reffered
to as the $f,g$-principal isotope of $(G,\oplus )$.

A subloop(subquasigroup) $(S,\otimes )$ of a loop(quasigroup)
$(G,\otimes )$ is called a Smarandache $f,g$-principal isotope of
the subloop(subquasigroup) $(S,\oplus )$ of a loop(quasigroup)
$(G,\oplus )$ if for some $f,g\in S$,
\begin{displaymath}
xR_g\otimes yL_f=(x\oplus y)~\forall~x,y\in S.
\end{displaymath}
On the other hand $(G,\otimes )$ is called a Smarandache
$f,g$-principal isotope of $(G,\oplus )$ if for some $f,g\in S$,
\begin{displaymath}
xR_g\otimes yL_f=(x\oplus y)~\forall~x,y\in G
\end{displaymath}
where $(S,\oplus )$ is a S-subquasigroup(S-subloop) of $(G,\oplus
)$. In these cases, $f$ and $g$ are called Smarandache
elements(S-elements).
\end{mydef}

\begin{myth}\label{1:1}(\cite{phd41})
Let $(G,\oplus)$ and $(H,\otimes)$ be two distinct isotopic
loops(quasigroups). There exists an $f,g$-principal isotope
$(G,\circ )$ of $(G,\oplus)$ such that $(H,\otimes)\cong (G,\circ
)$.
\end{myth}

\begin{mycor}\label{1:2}
Let ${\cal P}$ be an isotopic invariant property in
loops(quasigroups). If $(G,\oplus)$ is a loop(quasigroup) with the
property ${\cal P}$, then $(G,\oplus)$ is a universal
loop(quasigroup) relative to the property ${\cal P}$ if and only if
every $f,g$-principal isotope $(G,\circ )$ of $(G,\oplus)$ has the
property ${\cal P}$.
\end{mycor}
{\bf Proof}\\
If $(G,\oplus)$ is a universal loop relative to the property ${\cal
P}$ then every distinct loop isotope $(H,\otimes)$ of $(G,\oplus)$
has the property ${\cal P}$. By Theorem~\ref{1:1}, there exists an
$f,g$-principal isotope $(G,\circ )$ of $(G,\oplus)$ such that
$(H,\otimes)\cong (G,\circ )$. Hence, since ${\cal P}$ is an
isomorphic invariant property, every $(G,\circ )$ has it.\\

Conversely, if every $f,g$-principal isotope $(G,\circ )$ of
$(G,\oplus)$ has the property ${\cal P}$ and since by
Theorem~\ref{1:1} for each distinct isotope $(H,\otimes)$ there
exists an $f,g$-principal isotope $(G,\circ )$ of $(G,\oplus)$ such
that $(H,\otimes)\cong (G,\circ )$, then all $(H,\otimes)$ has the
property, Thus, $(G,\oplus)$ is a universal loop relative to the
property ${\cal P}$.

\begin{mylem}\label{1:3}
Let $(G,\oplus)$ be a loop(quasigroup) with a subloop(subquasigroup)
$(S,\oplus )$. If $(G,\circ )$ is an arbitrary $f,g$-principal
isotope of $(G,\oplus)$, then $(S,\circ )$ is a
subloop(subquasigroup) of $(G,\circ)$ if $(S,\circ )$ is a
Smarandache $f,g$-principal isotope of $(S,\oplus )$.
\end{mylem}
{\bf Proof}\\
If $(S,\circ )$ is a Smarandache $f,g$-principal isotope of
$(S,\oplus )$, then for some $f,g\in S$,
\begin{displaymath}
xR_g\circ yL_f=(x\oplus y)~\forall~x,y\in S\Rightarrow x\circ
y=xR_g^{-1}\oplus yL_f^{-1}\in S~\forall~x,y\in S
\end{displaymath}
since $f,g\in S$. So, $(S,\circ )$ is a subgroupoid of $(G,\circ )$.
$(S,\circ )$ is a subquasigroup follows from the fact that
$(S,\oplus )$ is a subquasigroup. $f\oplus g$ is a two sided
identity element in $(S,\circ )$. Thus, $(S,\circ )$ is a subloop of
$(G,\circ )$.

\section{Main Results}
\subsection*{Universality of Smarandache Loops}
\begin{myth}\label{1:4}
A Smarandache quasigroup is universal if all its $f,g$-principal
isotopes are Smarandache $f,g$-principal isotopes.
\end{myth}
{\bf Proof}\\
Let $(G,\oplus)$ be a Smarandache quasigroup with a S-subquasigroup
$(S,\oplus )$. If $(G,\circ )$ is an arbitrary $f,g$-principal
isotope of $(G,\oplus)$, then by Lemma~\ref{1:3}, $(S,\circ )$ is a
subquasigroup of $(G,\circ)$ if $(S,\circ )$ is a Smarandache
$f,g$-principal isotope of $(S,\oplus )$. Let us choose all
$(S,\circ )$ in this manner. So,
\begin{displaymath}
x\circ y=xR_g^{-1}\oplus yL_f^{-1}~\forall~x,y\in S.
\end{displaymath} It shall now be shown that
\begin{displaymath}(x\circ y)\circ z=x\circ (y\circ z)~\forall~x,y,z\in
S.
\end{displaymath}
But in the quasigroup $(G,\oplus )$, $xy$ will have preference over
$x\oplus y~\forall~x,y\in G$.
\begin{displaymath}
(x\circ y)\circ z=(xR_g^{-1}\oplus yL_f^{-1})\circ z=(xg^{-1}\oplus
f^{-1}y)\circ z=(xg^{-1}\oplus f^{-1}y)R_g^{-1}\oplus
zL_f^{-1}
\end{displaymath}
\begin{displaymath}
=(xg^{-1}\oplus f^{-1}y)g^{-1}\oplus f^{-1}z=xg^{-1}\oplus
f^{-1}yg^{-1}\oplus f^{-1}z.
\end{displaymath}
\begin{displaymath}
x\circ (y\circ z)=x\circ (yR_g^{-1}\oplus zL_f^{-1})=x\circ
(yg^{-1}\oplus f^{-1}z)=xR_g^{-1}\oplus (yg^{-1}\oplus
f^{-1}z)L_f^{-1}
\end{displaymath}
\begin{displaymath}
=xg^{-1}\oplus f^{-1}(yg^{-1}\oplus
f^{-1}z)=xg^{-1}\oplus f^{-1}yg^{-1}\oplus f^{-1}z.
\end{displaymath}

Thus, $(S,\circ )$ is an S-subquasigroup of $(G,\circ )$ hence,
$(G,\circ )$ is a S-quasigroup. By Theorem~\ref{1:1}, for any
isotope $(H,\otimes )$ of $(G,\oplus)$, there exists a $(G,\circ )$
such that $(H,\otimes )\cong (G,\circ )$. So we can now choose the
isomorphic image of $(S,\circ)$ which will now be an S-subquasigroup
in $(H,\otimes )$. So, $(H,\otimes )$ is an S-quasigroup. This
conclusion can also be drawn straight from Corollary~\ref{1:2}.

\begin{myth}\label{1:5}
A Smarandache loop is universal if all its $f,g$-principal isotopes
are Smarandache $f,g$-principal isotopes. Conversely, if a
Smarandache loop is universal then
\begin{displaymath}
(I,L_fR_g^{-1}R_{f^\rho}L_f^{-1},R_g^{-1}R_{f^\rho})
\end{displaymath} is an autotopism of an S-subloop of the S-loop such that $f$ and $g$ are S-elements.
\end{myth}
{\bf Proof}\\
Every loop is a quasigroup. Hence, the first claim follows from
Theorem~\ref{1:4}. The proof of the converse is as follows. If a
Smarandache loop $(G,\oplus )$ is universal then every isotope
$(H,\otimes)$ is an S-loop i.e there exists an S-subloop $(S,\otimes
)$ in $(H,\otimes )$. Let $(G,\circ )$ be the $f,g$-principal
isotope of $(G,\oplus)$, then by Corollary~\ref{1:2}, $(G,\circ)$ is
an S-loop with say an S-subloop $(S,\circ)$. So,
\begin{displaymath}
(x\circ y)\circ z=x\circ (y\circ z)~\forall~x,y,z\in S
\end{displaymath}
where \begin{displaymath} x\circ y=xR_g^{-1}\oplus
yL_f^{-1}~\forall~x,y\in S.
\end{displaymath}
\begin{displaymath}
(xR_g^{-1}\oplus yL_f^{-1})R_g^{-1}\oplus zL_f^{-1}=xR_g^{-1}\oplus
(yR_g^{-1}\oplus zL_f^{-1})L_f^{-1}. \end{displaymath} Replacing
$xR_g^{-1}$ by $x'$, $yL_f^{-1}$ by $y'$ and taking $z=e$ in
$(S,\oplus)$ we have; \begin{displaymath} (x'\oplus
y')R_g^{-1}R_{f^\rho}=x'\oplus
y'L_fR_g^{-1}R_{f^\rho}L_f^{-1}\Rightarrow (I,L_fR_
g^{-1}R_{f^\rho}L_f^{-1},R_g^{-1}R_{f^\rho}) \end{displaymath} is an
autotopism of an S-subloop $(S,\oplus )$ of the S-loop $(G,\oplus )$
such that $f$ and $g$ are S-elements.

\subsection*{Universality of Smarandache Bol, Moufang and Extra Loops}
\begin{myth}\label{1:6}
A Smarandache right(left)Bol loop is universal if all its
$f,g$-principal isotopes are Smarandache $f,g$-principal isotopes.
Conversely, if a Smarandache right(left)Bol loop is universal then
\begin{displaymath}
{\cal
T}_1=(R_gR_{f^\rho}^{-1},L_{g^\lambda}R_g^{-1}R_{f^\rho}L_f^{-1},R_g^{-1}R_{f^\rho})\bigg({\cal
T}_2=(R_{f^\rho}L_f^{-1}L_{g^\lambda}R_g^{-1},L_fL_{g^\lambda}^{-1},L_f^{-1}L_{g^\lambda})\bigg)
\end{displaymath}
is an autotopism of an SRB(SLB)-subloop of the SRBL(SLBL) such that
$f$ and $g$ are S-elements.
\end{myth}
{\bf Proof}\\
Let $(G,\oplus)$ be a SRBL(SLBL) with a S-RB(LB)-subloop $(S,\oplus
)$. If $(G,\circ )$ is an arbitrary $f,g$-principal isotope of
$(G,\oplus)$, then by Lemma~\ref{1:3}, $(S,\circ )$ is a subloop of
$(G,\circ)$ if $(S,\circ )$ is a Smarandache $f,g$-principal isotope
of $(S,\oplus )$. Let us choose all $(S,\circ )$ in this manner. So,
\begin{displaymath}
x\circ y=xR_g^{-1}\oplus yL_f^{-1}~\forall~x,y\in S.
\end{displaymath}
It is already known from \cite{phd3} that RB(LB) loops are
universal, hence $(S,\circ )$ is a RB(LB) loop thus an
S-RB(LB)-subloop of $(G,\circ)$. By Theorem~\ref{1:1}, for any
isotope $(H,\otimes )$ of $(G,\oplus)$, there exists a $(G,\circ )$
such that $(H,\otimes )\cong (G,\circ )$. So we can now choose the
isomorphic image of $(S,\circ)$ which will now be an
S-RB(LB)-subloop in $(H,\otimes )$. So, $(H,\otimes )$ is an
SRBL(SLBL). This conclusion can also be drawn straight from
Corollary~\ref{1:2}.

The proof of the converse is as follows. If a SRBL(SLBL) $(G,\oplus
)$ is universal then every isotope $(H,\otimes)$ is an SRBL(SLBL)
i.e there exists an S-RB(LB)-subloop $(S,\otimes )$ in $(H,\otimes
)$. Let $(G,\circ )$ be the $f,g$-principal isotope of $(G,\oplus)$,
then by Corollary~\ref{1:2}, $(G,\circ)$ is an SRBL(SLBL) with say
an SRB(SLB)-subloop $(S,\circ)$. So for an SRB-subloop $(S,\circ)$,
\begin{displaymath}
[(y\circ x)\circ z]\circ x=y\circ [(x\circ z)\circ
x]~\forall~x,y,z\in S\end{displaymath} where
\begin{displaymath}
x\circ y=xR_g^{-1}\oplus yL_f^{-1}~\forall~x,y\in S.
\end{displaymath}
Thus,
\begin{displaymath}
[(yR_g^{-1}\oplus xL_f^{-1})R_g^{-1}\oplus zL_f^{-1}]R_g^{-1}\oplus
xL_f^{-1}=yR_g^{-1}\oplus [(xR_g^{-1}\oplus zL_f^{-1})R_g^{-1}\oplus
xL_f^{-1}]L_f^{-1}. \end{displaymath} Replacing $yR_g^{-1}$ by $y'$,
$zL_f^{-1}$ by $z'$ and taking $x=e$ in $(S,\oplus)$ we have
\begin{displaymath}
(y'R_{f^\rho}R_g^{-1}\oplus z')R_g^{-1}R_{f^\rho}=y'\oplus
z'L_{g^\lambda}R_g^{-1}R_{f^\rho}L_f^{-1}. \end{displaymath} Again,
replace $y'R_{f^\rho}R_g^{-1}$ by $y''$ so that
\begin{displaymath}
(y''\oplus z')R_g^{-1}R_{f^\rho}=y''R_gR_{f^\rho}^{-1}\oplus
z'L_{g^\lambda}R_g^{-1}R_{f^\rho}L_f^{-1}\Rightarrow
(R_gR_{f^\rho}^{-1},L_{g^\lambda}R_g^{-1}R_{f^\rho}L_f^{-1},R_g^{-1}R_{f^\rho})
\end{displaymath}
is an autotopism of an SRB-subloop $(S,\oplus )$ of the S-loop $(G,\oplus )$ such that $f$ and $g$ are S-elements.\\
On the other hand, for a SLB-subloop $(S,\circ)$,
\begin{displaymath}
[x\circ (y\circ x)]\circ z=x\circ [y\circ (x\circ
z)]~\forall~x,y,z\in S\end{displaymath} where
\begin{displaymath}
x\circ y=xR_g^{-1}\oplus yL_f^{-1}~\forall~x,y\in S.
\end{displaymath}
Thus,
\begin{displaymath}
[xR_g^{-1}\oplus (yR_g^{-1}\oplus xL_f^{-1})L_f^{-1}]R_g^{-1}\oplus
zL_f^{-1}=xR_g^{-1}\oplus [yR_g^{-1}\oplus (xR_g^{-1}\oplus
zL_f^{-1})L_f^{-1}]L_f^{-1}.
\end{displaymath} Replacing $yR_g^{-1}$ by $y'$, $zL_f^{-1}$ by $z'$
and taking $x=e$ in $(S,\oplus)$ we have
\begin{displaymath}
y'R_{f^\rho}L_f^{-1}L_{g^\lambda}R_g^{-1}\oplus z'=(y'\oplus
z'L_{g^\lambda}L_f^{-1})L_f^{-1}L_{g^\lambda}.
\end{displaymath} Again, replace $z'L_{g^\lambda}L_f^{-1}$ by $z''$
so that
\begin{displaymath}
y'R_{f^\rho}L_f^{-1}L_{g^\lambda}R_g^{-1}\oplus
z''L_fL_{g^\lambda}^{-1}=(y'\oplus
z'')L_f^{-1}L_{g^\lambda}\Rightarrow
(R_{f^\rho}L_f^{-1}L_{g^\lambda}R_g^{-1},L_fL_{g^\lambda}^{-1},L_f^{-1}L_{g^\lambda})
\end{displaymath}
is an autotopism of an SLB-subloop $(S,\oplus )$ of the S-loop
$(G,\oplus )$ such that $f$ and $g$ are S-elements.

\begin{myth}\label{1:7}
A Smarandache Moufang loop is universal if all its $f,g$-principal
isotopes are Smarandache $f,g$-principal isotopes. Conversely, if a
Smarandache Moufang loop is universal then
\begin{displaymath}
(R_{g}L_f^{-1}L_{g^\lambda}R_g^{-1},L_fR_g^{-1}R_{f^\rho}L_f^{-1},L_f^{-1}L_{g^\lambda}R_g^{-1}R_{f^\rho}),
(R_{g}L_f^{-1}L_{g^\lambda}R_g^{-1},L_fR_g^{-1}R_{f^\rho}L_f^{-1},R_g^{-1}R_{f^\rho}L_f^{-1}L_{g^\lambda}),
\end{displaymath}
\begin{displaymath}
(R_gL_f^{-1}L_{g^\lambda}R_g^{-1}R_{f^\rho}R_g^{-1},L_fL_{g^\lambda}^{-1},L_f^{-1}L_{g^\lambda}),
(R_gR_{f^\rho}^{-1},L_fR_g^{-1}R_{f^\rho}L_f^{-1}L_{g^\lambda}L_f^{-1},R_g^{-1}R_{f^\rho}),
\end{displaymath}
\begin{displaymath}
(R_gL_f^{-1}L_{g^\lambda}R_g^{-1},L_{g^\lambda}R_g^{-1}R_{f^\rho}L_{g^\lambda}^{-1},R_g^{-1}R_{f^\rho}L_f^{-1}L_{g^\lambda}),
(R_{f^\rho}L_f^{-1}L_{g^\lambda}R_{f^\rho}^{-1},L_fR_g^{-1}R_{f^\rho}L_f^{-1},L_f^{-1}L_{g^\lambda}R_g^{-1}R_{f^\rho})
\end{displaymath}
are autotopisms of an SM-subloop of the SML such that $f$ and $g$
are S-elements.
\end{myth}
{\bf Proof}\\
Let $(G,\oplus)$ be a SML with a SM-subloop $(S,\oplus )$. If
$(G,\circ )$ is an arbitrary $f,g$-principal isotope of
$(G,\oplus)$, then by Lemma~\ref{1:3}, $(S,\circ )$ is a subloop of
$(G,\circ)$ if $(S,\circ )$ is a Smarandache $f,g$-principal isotope
of $(S,\oplus )$. Let us choose all $(S,\circ )$ in this manner. So,
\begin{displaymath}
x\circ y=xR_g^{-1}\oplus yL_f^{-1}~\forall~x,y\in S.
\end{displaymath}
It is already known from \cite{phd3} that Moufang loops are
universal, hence $(S,\circ )$ is a Moufang loop thus an SM-subloop
of $(G,\circ)$. By Theorem~\ref{1:1}, for any isotope $(H,\otimes )$
of $(G,\oplus)$, there exists a $(G,\circ )$ such that $(H,\otimes
)\cong (G,\circ )$. So we can now choose the isomorphic image of
$(S,\circ)$ which will now be an SM-subloop in $(H,\otimes )$. So,
$(H,\otimes )$ is an SML. This conclusion can also be drawn straight
from Corollary~\ref{1:2}.

The proof of the converse is as follows. If a SML $(G,\oplus )$ is
universal then every isotope $(H,\otimes)$ is an SML i.e there
exists an SM-subloop $(S,\otimes )$ in $(H,\otimes )$. Let $(G,\circ
)$ be the $f,g$-principal isotope of $(G,\oplus)$, then by
Corollary~\ref{1:2}, $(G,\circ)$ is an SML with say an SM-subloop
$(S,\circ)$. For an SM-subloop $(S,\circ)$,
\begin{displaymath}
(x\circ y)\circ (z\circ x)=[x\circ (y\circ z)]\circ
x~\forall~x,y,z\in S\end{displaymath} where
\begin{displaymath}
x\circ y=xR_g^{-1}\oplus yL_f^{-1}~\forall~x,y\in S.
\end{displaymath}
Thus,
\begin{displaymath}
(xR_g^{-1}\oplus yL_f^{-1})R_g^{-1}\oplus (zR_g^{-1}\oplus
xL_f^{-1})L_f^{-1}=[xR_g^{-1}\oplus (yR_g^{-1}\oplus
zL_f^{-1})L_f^{-1}]R_g^{-1}\oplus xL_f^{-1}.
\end{displaymath} Replacing $yR_g^{-1}$ by $y'$, $zL_f^{-1}$ by $z'$
and taking $x=e$ in $(S,\oplus)$ we have
\begin{displaymath}
y'R_gL_f^{-1}L_{g^\lambda}R_g^{-1}\oplus
z'L_fR_g^{-1}R_{f^\rho}L_f^{-1}=(y'\oplus
z')L_f^{-1}L_{g^\lambda}R_g^{-1}R_{f^\rho}\Rightarrow
\end{displaymath}
\begin{displaymath}
(R_{g}L_f^{-1}L_{g^\lambda}R_g^{-1},L_fR_g^{-1}R_{f^\rho}L_f^{-1},L_f^{-1}L_{g^\lambda}R_g^{-1}R_{f^\rho})
\end{displaymath}
is an autotopism of an SM-subloop $(S,\oplus )$ of the S-loop
$(G,\oplus )$ such that $f$ and $g$ are S-elements.

Again, for an SM-subloop $(S,\circ)$,
\begin{displaymath}
(x\circ y)\circ (z\circ x)=x\circ [(y\circ z)\circ x]
~\forall~x,y,z\in S\end{displaymath} where
\begin{displaymath}
x\circ y=xR_g^{-1}\oplus yL_f^{-1}~\forall~x,y\in S.
\end{displaymath}
Thus,
\begin{displaymath}
(xR_g^{-1}\oplus yL_f^{-1})R_g^{-1}\oplus (zR_g^{-1}\oplus
xL_f^{-1})L_f^{-1}=xR_g^{-1}\oplus [(yR_g^{-1}\oplus
zL_f^{-1})R_g^{-1}\oplus xL_f^{-1}]L_f^{-1}.
\end{displaymath} Replacing $yR_g^{-1}$ by $y'$, $zL_f^{-1}$ by $z'$
and taking $x=e$ in $(S,\oplus)$ we have
\begin{displaymath}
y'R_gL_f^{-1}L_{g^\lambda}R_g^{-1}\oplus
z'L_fR_g^{-1}R_{f^\rho}L_f^{-1}=(y'\oplus
z')R_g^{-1}R_{f^\rho}L_f^{-1}L_{g^\lambda}\Rightarrow
\end{displaymath}
\begin{displaymath}
(R_{g}L_f^{-1}L_{g^\lambda}R_g^{-1},L_fR_g^{-1}R_{f^\rho}L_f^{-1},R_g^{-1}R_{f^\rho}L_f^{-1}L_{g^\lambda})
\end{displaymath}
is an autotopism of an SM-subloop $(S,\oplus )$ of the S-loop
$(G,\oplus )$ such that $f$ and $g$ are S-elements.

Also, if $(S,\circ)$ is an SM-subloop then,
\begin{displaymath}
[(x\circ y)\circ x]\circ z=x\circ [y\circ (x\circ z)]
~\forall~x,y,z\in S\end{displaymath} where
\begin{displaymath}
x\circ y=xR_g^{-1}\oplus yL_f^{-1}~\forall~x,y\in S.
\end{displaymath}
Thus,
\begin{displaymath}
[(xR_g^{-1}\oplus yL_f^{-1})R_g^{-1}\oplus xL_f^{-1}]R_g^{-1}\oplus
zL_f^{-1}=xR_g^{-1}\oplus [yR_g^{-1}\oplus (xR_g^{-1}\oplus
zL_f^{-1})L_f^{-1}]L_f^{-1}.
\end{displaymath} Replacing $yR_g^{-1}$ by $y'$, $zL_f^{-1}$ by $z'$
and taking $x=e$ in $(S,\oplus)$ we have
\begin{displaymath}
y'R_gL_f^{-1}L_{g^\lambda}R_g^{-1}R_{f^\rho}R_g^{-1}\oplus
z'=(y'\oplus z'L_{g^\lambda}L_f^{-1})L_f^{-1}L_{g^\lambda}.
\end{displaymath}
Again, replace $z'L_{g^\lambda}L_f^{-1}$ by $z''$ so that
\begin{displaymath}
y'R_gL_f^{-1}L_{g^\lambda}R_g^{-1}R_{f^\rho}R_g^{-1}\oplus
z''L_fL_{g^\lambda}^{-1}=(y'\oplus
z'')L_f^{-1}L_{g^\lambda}\Rightarrow
(R_gL_f^{-1}L_{g^\lambda}R_g^{-1}R_{f^\rho}R_g^{-1},L_fL_{g^\lambda}^{-1},L_f^{-1}L_{g^\lambda})
\end{displaymath}
is an autotopism of an SM-subloop $(S,\oplus )$ of the S-loop
$(G,\oplus )$ such that $f$ and $g$ are S-elements.

Furthermore, if $(S,\circ)$ is an SM-subloop then,
\begin{displaymath}
[(y\circ x)\circ z]\circ x=y\circ [x\circ (z\circ x)]
~\forall~x,y,z\in S\end{displaymath} where
\begin{displaymath}
x\circ y=xR_g^{-1}\oplus yL_f^{-1}~\forall~x,y\in S.
\end{displaymath}
Thus,
\begin{displaymath}
[(yR_g^{-1}\oplus xL_f^{-1})R_g^{-1}\oplus zL_f^{-1}]R_g^{-1}\oplus
xL_f^{-1}=yR_g^{-1}\oplus [xR_g^{-1}\oplus (zR_g^{-1}\oplus
xL_f^{-1})L_f^{-1}]L_f^{-1}.
\end{displaymath} Replacing $yR_g^{-1}$ by $y'$, $zL_f^{-1}$ by $z'$
and taking $x=e$ in $(S,\oplus)$ we have
\begin{displaymath}
(y'R_{f^\rho}R_g^{-1}\oplus z')R_g^{-1}R_{f^\rho}=y'\oplus
z'L_fR_g^{-1}R_{f^\rho}L_f^{-1}L_{g^\lambda}L_f^{-1}.
\end{displaymath}
Again, replace $y'R_{f^\rho}R_g^{-1}$ by $y''$ so that
\begin{displaymath}
(y''\oplus z')R_g^{-1}R_{f^\rho}=y''R_gR_{f^\rho}^{-1}\oplus
z'L_fR_g^{-1}R_{f^\rho}L_f^{-1}L_{g^\lambda}L_f^{-1}\Rightarrow
(R_gR_{f^\rho}^{-1},L_fR_g^{-1}R_{f^\rho}L_f^{-1}L_{g^\lambda}L_f^{-1},R_g^{-1}R_{f^\rho})
\end{displaymath}
is an autotopism of an SM-subloop $(S,\oplus )$ of the S-loop
$(G,\oplus )$ such that $f$ and $g$ are S-elements.

Lastly, $(S,\oplus)$ is an SM-subloop if and only if $(S,\circ)$ is
an SRB-subloop and an SLB-subloop. So by Theorem~\ref{1:6}, ${\cal
T}_1$ and ${\cal T}_2$ are autotopisms in $(S,\oplus)$, hence ${\cal
T}_1{\cal T}_2$ and ${\cal T}_2{\cal T}_1$ are autotopisms in
$(S,\oplus)$.

\begin{myth}\label{1:8}
A Smarandache extra loop is universal if all its $f,g$-principal
isotopes are Smarandache $f,g$-principal isotopes. Conversely, if a
Smarandache extra loop is universal then
$(R_gL_f^{-1}L_{g^\lambda}R_g^{-1},L_fR_{f^\rho}^{-1}R_gL_f^{-1},L_f^{-1}L_{g^\lambda}R_{f^\rho}^{-1}R_g)$,
\begin{displaymath}
(R_gR_{f^\rho}^{-1}R_gL_f^{-1}L_{g^\lambda}R_g^{-1},L_{g^\lambda}L_f^{-1},L_f^{-1}L_{g^\lambda}),
(R_{f^\rho}R_g^{-1},L_fL_{g^\lambda}^{-1}L_fR_g^{-1}R_{f^\rho}L_f^{-1},R_g^{-1}R_{f^\rho})
\end{displaymath}
\begin{displaymath}
(R_{g}L_f^{-1}L_{g^\lambda}R_g^{-1},L_fR_g^{-1}R_{f^\rho}L_f^{-1},L_f^{-1}L_{g^\lambda}R_g^{-1}R_{f^\rho}),
(R_{g}L_f^{-1}L_{g^\lambda}R_g^{-1},L_fR_g^{-1}R_{f^\rho}L_f^{-1},R_g^{-1}R_{f^\rho}L_f^{-1}L_{g^\lambda}),
\end{displaymath}
\begin{displaymath}
(R_gL_f^{-1}L_{g^\lambda}R_g^{-1}R_{f^\rho}R_g^{-1},L_fL_{g^\lambda}^{-1},L_f^{-1}L_{g^\lambda}),
(R_gR_{f^\rho}^{-1},L_fR_g^{-1}R_{f^\rho}L_f^{-1}L_{g^\lambda}L_f^{-1},R_g^{-1}R_{f^\rho}),
\end{displaymath}
\begin{displaymath}
(R_gL_f^{-1}L_{g^\lambda}R_g^{-1},L_{g^\lambda}R_g^{-1}R_{f^\rho}L_{g^\lambda}^{-1},R_g^{-1}R_{f^\rho}L_f^{-1}L_{g^\lambda}),
(R_{f^\rho}L_f^{-1}L_{g^\lambda}R_{f^\rho}^{-1},L_fR_g^{-1}R_{f^\rho}L_f^{-1},L_f^{-1}L_{g^\lambda}R_g^{-1}R_{f^\rho}),
\end{displaymath}
are autotopisms of an SE-subloop of the SEL such that $f$ and $g$
are S-elements.
\end{myth}
{\bf Proof}\\
Let $(G,\oplus)$ be a SEL with a SE-subloop $(S,\oplus )$. If
$(G,\circ )$ is an arbitrary $f,g$-principal isotope of
$(G,\oplus)$, then by Lemma~\ref{1:3}, $(S,\circ )$ is a subloop of
$(G,\circ)$ if $(S,\circ )$ is a Smarandache $f,g$-principal isotope
of $(S,\oplus )$. Let us choose all $(S,\circ )$ in this manner. So,
\begin{displaymath}
x\circ y=xR_g^{-1}\oplus yL_f^{-1}~\forall~x,y\in S.
\end{displaymath}
In \cite{phd34} and \cite{phd36} respectively, it was shown and
stated that a loop is an extra loop if and only if it is a Moufang
loop and a CC-loop. But since CC-loops are G-loops(they are
isomorphic to all loop isotopes) then extra loops are universal,
hence $(S,\circ )$ is an extra loop thus an SE-subloop of
$(G,\circ)$. By Theorem~\ref{1:1}, for any isotope $(H,\otimes )$ of
$(G,\oplus)$, there exists a $(G,\circ )$ such that $(H,\otimes
)\cong (G,\circ )$. So we can now choose the isomorphic image of
$(S,\circ)$ which will now be an SE-subloop in $(H,\otimes )$. So,
$(H,\otimes )$ is an SEL. This conclusion can also be drawn straight
from Corollary~\ref{1:2}.

The proof of the converse is as follows. If a SEL $(G,\oplus )$ is
universal then every isotope $(H,\otimes)$ is an SEL i.e there
exists an SE-subloop $(S,\otimes )$ in $(H,\otimes )$. Let $(G,\circ
)$ be the $f,g$-principal isotope of $(G,\oplus)$, then by
Corollary~\ref{1:2}, $(G,\circ)$ is an SEL with say an SE-subloop
$(S,\circ)$. For an SE-subloop $(S,\circ)$,
\begin{displaymath}
[(x\circ y)\circ z]\circ x=x\circ [y\circ (z\circ
x)]~\forall~x,y,z\in S\end{displaymath} where
\begin{displaymath}
x\circ y=xR_g^{-1}\oplus yL_f^{-1}~\forall~x,y\in S.
\end{displaymath}
Thus,
\begin{displaymath}
[(xR_g^{-1}\oplus yL_f^{-1})R_g^{-1}\oplus zL_f^{-1}]R_g^{-1}\oplus
xL_f^{-1}=xR_g^{-1}\oplus [yR_g^{-1}\oplus (zR_g^{-1}\oplus
xL_f^{-1})L_f^{-1}]L_f^{-1}.
\end{displaymath} Replacing $yR_g^{-1}$ by $y'$, $zL_f^{-1}$ by $z'$
and taking $x=e$ in $(S,\oplus)$ we have
\begin{displaymath}
(y'R_gL_f^{-1}L_{g^\lambda}R_g^{-1}\oplus
z')R_g^{-1}R_{f^\rho}=(y'\oplus
z'L_fR_g^{-1}R_{f^\rho}L_f^{-1})L_f^{-1}L_{g^\lambda}.
\end{displaymath}
Again, replace $z'L_fR_g^{-1}R_{f^\rho}L_f^{-1}$ by $z''$ so that
\begin{displaymath}
y'R_gL_f^{-1}L_{g^\lambda}R_g^{-1}\oplus
z''L_fR_{f^\rho}^{-1}R_gL_f^{-1}=(y'\oplus
z'')L_f^{-1}L_{g^\lambda}R_{f^\rho}^{-1}R_g\Rightarrow
\end{displaymath}
\begin{displaymath}
(R_gL_f^{-1}L_{g^\lambda}R_g^{-1},L_fR_{f^\rho}^{-1}R_gL_f^{-1},L_f^{-1}L_{g^\lambda}R_{f^\rho}^{-1}R_g)
\end{displaymath}
is an autotopism of an SE-subloop $(S,\oplus )$ of the S-loop
$(G,\oplus )$ such that $f$ and $g$ are S-elements.

Again, for an SE-subloop $(S,\circ)$,
\begin{displaymath}
(x\circ y)\circ (x\circ z)=x\circ [(y\circ x)\circ z]
~\forall~x,y,z\in S\end{displaymath} where
\begin{displaymath}
x\circ y=xR_g^{-1}\oplus yL_f^{-1}~\forall~x,y\in S.
\end{displaymath}
Thus,
\begin{displaymath}
(xR_g^{-1}\oplus yL_f^{-1})R_g^{-1}\oplus (xR_g^{-1}\oplus
zL_f^{-1})L_f^{-1}=xR_g^{-1}\oplus [(yR_g^{-1}\oplus
xL_f^{-1})R_g^{-1}\oplus zL_f^{-1}]L_f^{-1}.
\end{displaymath} Replacing $yR_g^{-1}$ by $y'$, $zL_f^{-1}$ by $z'$
and taking $x=e$ in $(S,\oplus)$ we have
\begin{displaymath}
y'R_gL_f^{-1}L_{g^\lambda}R_g^{-1}\oplus
z'L_{g^\lambda}L_f^{-1}=(y'R_{f^\rho}R_g^{-1}\oplus
z')L_f^{-1}L_{g^\lambda}.
\end{displaymath}
Again, replace $y'R_{f^\rho}R_g^{-1}$ by $y''$ so that
\begin{displaymath}
y''R_gR_{f^\rho}^{-1}R_gL_f^{-1}L_{g^\lambda}R_g^{-1}\oplus
z'L_{g^\lambda}L_f^{-1}=(y''\oplus
z')L_f^{-1}L_{g^\lambda}\Rightarrow
(R_gR_{f^\rho}^{-1}R_gL_f^{-1}L_{g^\lambda}R_g^{-1},L_{g^\lambda}L_f^{-1},L_f^{-1}L_{g^\lambda})
\end{displaymath}
is an autotopism of an SE-subloop $(S,\oplus )$ of the S-loop
$(G,\oplus )$ such that $f$ and $g$ are S-elements.

Also, if $(S,\circ)$ is an SE-subloop then,
\begin{displaymath}
(y\circ x)\circ (z\circ x)=[y\circ (x\circ z)]\circ x
~\forall~x,y,z\in S\end{displaymath} where
\begin{displaymath}
x\circ y=xR_g^{-1}\oplus yL_f^{-1}~\forall~x,y\in S.
\end{displaymath}
Thus,
\begin{displaymath}
(yR_g^{-1}\oplus xL_f^{-1})R_g^{-1}\oplus (zR_g^{-1}\oplus
xL_f^{-1})L_f^{-1}= [(yR_g^{-1}\oplus (xR_g^{-1}\oplus
zL_f^{-1})L_f^{-1}]R_g^{-1}\oplus xL_f^{-1}.
\end{displaymath} Replacing $yR_g^{-1}$ by $y'$, $zL_f^{-1}$ by $z'$
and taking $x=e$ in $(S,\oplus)$ we have
\begin{displaymath}
y'R_{f^\rho}R_g^{-1}\oplus z'L_fR_g^{-1}R_{f^\rho}L_f^{-1}=(y'\oplus
z'L_{g^\lambda}L_f^{-1})R_g^{-1}R_{f^\rho}.
\end{displaymath}
Again, replace $z'L_{g^\lambda}L_f^{-1}$ by $z''$ so that
\begin{displaymath}
y'R_{f^\rho}R_g^{-1}\oplus
z''L_fL_{g^\lambda}^{-1}L_fR_g^{-1}R_{f^\rho}L_f^{-1}=(y'\oplus
z')R_g^{-1}R_{f^\rho}\Rightarrow
(R_{f^\rho}R_g^{-1},L_fL_{g^\lambda}^{-1}L_fR_g^{-1}R_{f^\rho}L_f^{-1},R_g^{-1}R_{f^\rho})
\end{displaymath}
is an autotopism of an SE-subloop $(S,\oplus )$ of the S-loop
$(G,\oplus )$ such that $f$ and $g$ are S-elements.

Lastly, $(S,\oplus)$ is an SE-subloop if and only if $(S,\circ)$ is
an SM-subloop and an SCC-subloop. So by Theorem~\ref{1:7}, the six
remaining triples are autotopisms in $(S,\oplus)$.

\subsection*{Universality of Smarandache Inverse Property Loops}
\begin{myth}\label{1:9}
A Smarandache left(right) inverse property loop in which all its
$f,g$-principal isotopes are Smarandache $f,g$-principal isotopes is
universal if and only if it is a Smarandache left(right) Bol loop in
which all its $f,g$-principal isotopes are Smarandache
$f,g$-principal isotopes.
\end{myth}
{\bf Proof}\\
Let $(G,\oplus)$ be a SLIPL with a SLIP-subloop $(S,\oplus )$. If
$(G,\circ )$ is an arbitrary $f,g$-principal isotope of
$(G,\oplus)$, then by Lemma~\ref{1:3}, $(S,\circ )$ is a subloop of
$(G,\circ)$ if $(S,\circ )$ is a Smarandache $f,g$-principal isotope
of $(S,\oplus )$. Let us choose all $(S,\circ )$ in this manner. So,
\begin{displaymath}
x\circ y=xR_g^{-1}\oplus yL_f^{-1}~\forall~x,y\in S.
\end{displaymath}
$(G,\oplus)$ is a universal SLIPL if and only if every isotope
$(H,\otimes )$ is a SLIPL. $(H,\otimes )$ is a SLIPL if and only if
it has at least a SLIP-subloop $(S,\otimes )$. By Theorem~\ref{1:1},
for any isotope $(H,\otimes )$ of $(G,\oplus)$, there exists a
$(G,\circ )$ such that $(H,\otimes )\cong (G,\circ )$. So we can now
choose the isomorphic image of $(S,\circ)$ to be $(S,\otimes )$
which is already a SLIP-subloop in $(H,\otimes )$. So, $(S,\circ)$
is also a SLIP-subloop in $(G,\circ )$. As shown in \cite{phd3},
$(S,\oplus )$ and its $f,g$-isotope(Smarandache $f,g$-isotope)
$(S,\circ)$ are SLIP-subloops if and only if $(S,\oplus )$ is a left
Bol subloop(i.e a SLB-subloop). So, $(G,\oplus)$ is  SLBL.

Conversely, if  $(G,\oplus)$ is  SLBL, then there exists a
SLB-subloop $(S,\oplus )$ in $(G,\oplus)$. If $(G,\circ )$ is an
arbitrary $f,g$-principal isotope of $(G,\oplus)$, then by
Lemma~\ref{1:3}, $(S,\circ )$ is a subloop of $(G,\circ)$ if
$(S,\circ )$ is a Smarandache $f,g$-principal isotope of $(S,\oplus
)$. Let us choose all $(S,\circ )$ in this manner. So,
\begin{displaymath}
x\circ y=xR_g^{-1}\oplus yL_f^{-1}~\forall~x,y\in S.
\end{displaymath}
By Theorem~\ref{1:1}, for any isotope $(H,\otimes )$ of
$(G,\oplus)$, there exists a $(G,\circ )$ such that $(H,\otimes
)\cong (G,\circ )$. So we can now choose the isomorphic image of
$(S,\circ)$ to be $(S,\otimes )$ which is a SLB-subloop in
$(H,\otimes )$ using the same reasoning in Theorem~\ref{1:6}. So,
$(S,\circ)$ is a SLB-subloop in $(G,\circ )$. Left Bol loops have
the left inverse property(LIP), hence, $(S,\oplus )$ and $(S,\circ)$
are SLIP-subloops in $(G,\oplus)$ and $(G,\circ )$ respectively.
Thence, $(G,\oplus)$ and $(G,\circ )$ are SLBLs. Therefore,
$(G,\oplus)$ is
a universal SLIPL.\\

The proof for a Smarandache right inverse property loop is similar
and is as follows. Let $(G,\oplus)$ be a SRIPL with a SRIP-subloop
$(S,\oplus )$. If $(G,\circ )$ is an arbitrary $f,g$-principal
isotope of $(G,\oplus)$, then by Lemma~\ref{1:3}, $(S,\circ )$ is a
subloop of $(G,\circ)$ if $(S,\circ )$ is a Smarandache
$f,g$-principal isotope of $(S,\oplus )$. Let us choose all
$(S,\circ )$ in this manner. So,
\begin{displaymath}
x\circ y=xR_g^{-1}\oplus yL_f^{-1}~\forall~x,y\in S.
\end{displaymath}
$(G,\oplus)$ is a universal SRIPL if and only if every isotope
$(H,\otimes )$ is a SRIPL. $(H,\otimes )$ is a SRIPL if and only if
it has at least a SRIP-subloop $(S,\otimes )$. By Theorem~\ref{1:1},
for any isotope $(H,\otimes )$ of $(G,\oplus)$, there exists a
$(G,\circ )$ such that $(H,\otimes )\cong (G,\circ )$. So we can now
choose the isomorphic image of $(S,\circ)$ to be $(S,\otimes )$
which is already a SRIP-subloop in $(H,\otimes )$. So, $(S,\circ)$
is also a SRIP-subloop in $(G,\circ )$. As shown in \cite{phd3},
$(S,\oplus )$ and its $f,g$-isotope(Smarandache $f,g$-isotope)
$(S,\circ)$ are SRIP-subloops if and only if $(S,\oplus )$ is a
right Bol subloop(i.e a SRB-subloop). So, $(G,\oplus)$ is  SRBL.

Conversely, if  $(G,\oplus)$ is  SRBL, then there exists a
SRB-subloop $(S,\oplus )$ in $(G,\oplus)$. If $(G,\circ )$ is an
arbitrary $f,g$-principal isotope of $(G,\oplus)$, then by
Lemma~\ref{1:3}, $(S,\circ )$ is a subloop of $(G,\circ)$ if
$(S,\circ )$ is a Smarandache $f,g$-principal isotope of $(S,\oplus
)$. Let us choose all $(S,\circ )$ in this manner. So,
\begin{displaymath}
x\circ y=xR_g^{-1}\oplus yL_f^{-1}~\forall~x,y\in S.
\end{displaymath}
By Theorem~\ref{1:1}, for any isotope $(H,\otimes )$ of
$(G,\oplus)$, there exists a $(G,\circ )$ such that $(H,\otimes
)\cong (G,\circ )$. So we can now choose the isomorphic image of
$(S,\circ)$ to be $(S,\otimes )$ which is a SRB-subloop in
$(H,\otimes )$ using the same reasoning in Theorem~\ref{1:6}. So,
$(S,\circ)$ is a SRB-subloop in $(G,\circ )$. Right Bol loops have
the right inverse property(RIP), hence, $(S,\oplus )$ and
$(S,\circ)$ are SRIP-subloops in $(G,\oplus)$ and $(G,\circ )$
respectively. Thence, $(G,\oplus)$ and $(G,\circ )$ are SRBLs.
Therefore, $(G,\oplus)$ is a universal SRIPL.

\begin{myth}\label{1:10}
A Smarandache inverse property loop in which all its $f,g$-principal
isotopes are Smarandache $f,g$-principal isotopes is universal if
and only if it is a Smarandache Moufang loop in which all its
$f,g$-principal isotopes are Smarandache $f,g$-principal isotopes.
\end{myth}
{\bf Proof}\\
Let $(G,\oplus)$ be a SIPL with a SIP-subloop $(S,\oplus )$. If
$(G,\circ )$ is an arbitrary $f,g$-principal isotope of
$(G,\oplus)$, then by Lemma~\ref{1:3}, $(S,\circ )$ is a subloop of
$(G,\circ)$ if $(S,\circ )$ is a Smarandache $f,g$-principal isotope
of $(S,\oplus )$. Let us choose all $(S,\circ )$ in this manner. So,
\begin{displaymath}
x\circ y=xR_g^{-1}\oplus yL_f^{-1}~\forall~x,y\in S.
\end{displaymath}
$(G,\oplus)$ is a universal SIPL if and only if every isotope
$(H,\otimes )$ is a SIPL. $(H,\otimes )$ is a SIPL if and only if it
has at least a SIP-subloop $(S,\otimes )$. By Theorem~\ref{1:1}, for
any isotope $(H,\otimes )$ of $(G,\oplus)$, there exists a $(G,\circ
)$ such that $(H,\otimes )\cong (G,\circ )$. So we can now choose
the isomorphic image of $(S,\circ)$ to be $(S,\otimes )$ which is
already a SIP-subloop in $(H,\otimes )$. So, $(S,\circ)$ is also a
SIP-subloop in $(G,\circ )$. As shown in \cite{phd3}, $(S,\oplus )$
and its $f,g$-isotope(Smarandache $f,g$-isotope) $(S,\circ)$ are
SIP-subloops if and only if $(S,\oplus )$ is a Moufang subloop(i.e a
SM-subloop). So, $(G,\oplus)$ is  SML.

Conversely, if  $(G,\oplus)$ is SML, then there exists a SM-subloop
$(S,\oplus )$ in $(G,\oplus)$. If $(G,\circ )$ is an arbitrary
$f,g$-principal isotope of $(G,\oplus)$, then by Lemma~\ref{1:3},
$(S,\circ )$ is a subloop of $(G,\circ)$ if $(S,\circ )$ is a
Smarandache $f,g$-principal isotope of $(S,\oplus )$. Let us choose
all $(S,\circ )$ in this manner. So,
\begin{displaymath}
x\circ y=xR_g^{-1}\oplus yL_f^{-1}~\forall~x,y\in S.
\end{displaymath}
By Theorem~\ref{1:1}, for any isotope $(H,\otimes )$ of
$(G,\oplus)$, there exists a $(G,\circ )$ such that $(H,\otimes
)\cong (G,\circ )$. So we can now choose the isomorphic image of
$(S,\circ)$ to be $(S,\otimes )$ which is a SM-subloop in
$(H,\otimes )$ using the same reasoning in Theorem~\ref{1:6}. So,
$(S,\circ)$ is a SM-subloop in $(G,\circ )$. Moufang loops have the
inverse property(IP), hence, $(S,\oplus )$ and $(S,\circ)$ are
SIP-subloops in $(G,\oplus)$ and $(G,\circ )$ respectively. Thence,
$(G,\oplus)$ and $(G,\circ )$ are SMLs. Therefore, $(G,\oplus)$ is a
universal SIPL.

\begin{mycor}\label{1:11}
If a Smarandache left(right) inverse property loop is universal then
\begin{displaymath}
(R_gR_{f^\rho}^{-1},L_{g^\lambda}R_g^{-1}R_{f^\rho}L_f^{-1},R_g^{-1}R_{f^\rho})\bigg(
(R_{f^\rho}L_f^{-1}L_{g^\lambda}R_g^{-1},L_fL_{g^\lambda}^{-1},L_f^{-1}L_{g^\lambda})\bigg)
\end{displaymath}
is an autotopism of an SLIP(SRIP)-subloop of the SLIPL(SRIPL) such
that $f$ and $g$ are S-elements.
\end{mycor}
{\bf Proof}\\
This follows by Theorem~\ref{1:9} and Theorem~\ref{1:11}.

\begin{mycor}\label{1:12}
If a Smarandache inverse property loop is universal then
\begin{displaymath}
(R_{g}L_f^{-1}L_{g^\lambda}R_g^{-1},L_fR_g^{-1}R_{f^\rho}L_f^{-1},L_f^{-1}L_{g^\lambda}R_g^{-1}R_{f^\rho}),
(R_{g}L_f^{-1}L_{g^\lambda}R_g^{-1},L_fR_g^{-1}R_{f^\rho}L_f^{-1},R_g^{-1}R_{f^\rho}L_f^{-1}L_{g^\lambda}),
\end{displaymath}
\begin{displaymath}
(R_gL_f^{-1}L_{g^\lambda}R_g^{-1}R_{f^\rho}R_g^{-1},L_fL_{g^\lambda}^{-1},L_f^{-1}L_{g^\lambda}),
(R_gR_{f^\rho}^{-1},L_fR_g^{-1}R_{f^\rho}L_f^{-1}L_{g^\lambda}L_f^{-1},R_g^{-1}R_{f^\rho}),
\end{displaymath}
\begin{displaymath}
(R_gL_f^{-1}L_{g^\lambda}R_g^{-1},L_{g^\lambda}R_g^{-1}R_{f^\rho}L_{g^\lambda}^{-1},R_g^{-1}R_{f^\rho}L_f^{-1}L_{g^\lambda}),
(R_{f^\rho}L_f^{-1}L_{g^\lambda}R_{f^\rho}^{-1},L_fR_g^{-1}R_{f^\rho}L_f^{-1},L_f^{-1}L_{g^\lambda}R_g^{-1}R_{f^\rho})
\end{displaymath}
are autotopisms of an SIP-subloop of the SIPL such that $f$ and $g$
are S-elements.
\end{mycor}
{\bf Proof}\\
This follows from Theorem~\ref{1:10} and Theorem~\ref{1:7}.


\begin{thebibliography}{99}
\bibitem{phd30} R. Artzy (1959), {\it Crossed inverse and related
loops}, Trans. Amer. Math. Soc. 91, 3, 480--492.
\bibitem{phd41} R. H. Bruck (1966), {\it A survey of binary systems}, Springer-Verlag, Berlin-G\"ottingen-Heidelberg, 185pp.
\bibitem{phd40} R. H. Bruck and L. J. Paige (1956), {\it Loops whose
inner mappings are automorphisms}, The annuals of Mathematics, 63,
2, 308--323.
\bibitem{phd39} O. Chein, H. O. Pflugfelder and J. D. H. Smith (1990), {\it Quasigroups and loops : Theory and applications}, Heldermann Verlag, 568pp.
\bibitem{phd49} J. D\'{e}ne and A. D. Keedwell (1974), {\it Latin squares and their applications}, the Academic press Lts, 549pp.
\bibitem{phd50} F. Fenyves (1968), {\it Extra loops I}, Publ. Math. Debrecen, 15, 235--238.
\bibitem{phd56} F. Fenyves (1969), {\it Extra loops II}, Publ. Math. Debrecen, 16, 187--192.
\bibitem{phd42} E. G. Goodaire, E. Jespers and C. P. Milies (1996), {\it Alternative loop rings}, NHMS(184), Elsevier, 387pp.
\bibitem{phd34} E. G. Goodaire and D. A. Robinson (1990), {\it Some special conjugacy closed loops}, Canad. Math. Bull. 33, 73--78.
\bibitem{sma1} T. G. Ja\'iy\'e\d ol\'a (2006), {\it An holomorphic study of the Smarandache concept in
loops}, Scientia Magna Journal, 2, 1, 1--8.
\bibitem{sma2} T. G. Ja\'iy\'e\d ol\'a (2006), {\it Smarandache quasigroups}, Scientia Magna Journal, 2, 2, to appear.
\bibitem{phdtope} T. G. Ja\'iy\'e\d ol\'a (2005), {\it An isotopic study of
properties of central loops}, M.Sc. Dissertation, University of
Agriculture, Abeokuta.
\bibitem{tope} T. G. Ja\'iy\'e\d ol\'a and J. O. Ad\'en\'iran, {\it On isotopic characterization of central
loops}, communicated for publication.
\bibitem{phd118} T. Kepka, M. K. Kinyon, J. D. Phillips, {\it
F-quasigroups and generalised modules}, communicated for
publication.
\bibitem{phd119} T. Kepka, M. K. Kinyon, J. D. Phillips, {\it F-quasigroups isotopic to groups}, communicated for publication.
\bibitem{phd95} T.
Kepka, M. K. Kinyon, J. D. Phillips, {\it The structure of
F-quasigroups}, communicated for publication.
\bibitem{phd36} M. K. Kinyon, K. Kunen (2004), {\it The structure of
extra loops}, Quasigroups and Related Systems 12, 39--60.
\bibitem{phd124}  M. K. Kinyon, J. D. Phillips and P. Vojt\v echovsk\'y (2004), {\it Loops of Bol-Moufang type with a subgroup of index
two}, Bull. Acad. Stinte. Rebub. Mold. Mat. 2,45, 1--17.
\bibitem{ken2} K. Kunen (1996), {\it Quasigroups, loops and associative laws}, J. Alg.185, 194--204.
\bibitem{ken1} K. Kunen (1996), {\it Moufang quasigroups}, J. Alg. 183, 231--234.
\bibitem{muk} A. S. Muktibodh (2006), {\it Smarandache quasigroups},
Scientia Magna Journal, 2, 1, 13--19.
\bibitem{phd88} P. T. Nagy and K. Strambach (1994), {\it Loops as
invariant sections in groups, and their geometry}, Canad. J. Math.
46, 5, 1027--1056.
\bibitem{phd43} J. M. Osborn (1961), {\it Loops with the weak
inverse property}, Pac. J. Math. 10, 295--304.
\bibitem{phd3} H. O. Pflugfelder (1990), {\it Quasigroups and loops : Introduction}, Sigma series in Pure Math. 7, Heldermann Verlag, Berlin, 147pp.
\bibitem{phd9} J. D. Phillips and P. Vojt\v echovsk\'y (2005), {\it The varieties of loops of Bol-Moufang type}, Alg. Univer. (to appear).
\bibitem{phd61} J. D. Phillips and P. Vojt\v echovsk\'y (2005), {\it The varieties of quasigroups of Bol-Moufang type : An equational
approach}, J. Alg. 293, 17--33.
\bibitem{phd75} W. B. Vasantha Kandasamy (2002), {\it Smarandache
loops}, Department of Mathematics, Indian Institute of Technology,
Madras, India, 128pp.
\bibitem{phd83} W. B. Vasantha Kandasamy (2002), {\it Smarandache
Loops}, Smarandache notions journal, 13, 252--258.
\end{thebibliography}
\end{document}